\documentclass[fleqn]{mat01}
\usepackage{times,mathtimy,amssymb}
\begin{document}

\setcounter{page}{195}
\firstpage{195}

\font\xx=msam5 at 9pt
\def\ab{\mbox{\xx{\char'03}}}

\font\sa=tibi at 10.4pt

\def\defi{\trivlist\item[\hskip\labelsep{\bf DEFINITION.}]}
\def\remark{\trivlist\item[\hskip\labelsep{\it Remark.}]}
\def\remarks{\trivlist\item[\hskip\labelsep{\it Remarks}]}
\def\noot{\trivlist\item[\hskip\labelsep{{\it Note.}}]}
\def\thoe{\trivlist\item[\hskip\labelsep{{\bf Theorem}}]}

\newtheorem{dets}{\rm DEFINITION}
\newtheorem{theo}{Theorem}
\renewcommand\thetheo{\arabic{section}.\arabic{theo}}
\newtheorem{theor}[theo]{\bf Theorem}

\newcommand{\eps}{\varepsilon}
\newcommand{\To}{\longrightarrow}
\newcommand{\C}{\mathcal{C}}
\newcommand{\K}{\mathcal{K}}
\newcommand{\T}{\mathcal{T}}
\newcommand{\bq}{\begin{equation}}
\newcommand{\eq}{\end{equation}}

\def\pd#1#2{\frac{\partial #1}{\partial#2}}
\def\arrow{\to}

\title{Positive solutions of singular boundary value problem of~negative exponent
Emden--Fowler equation}

\markboth{Yuxia Wang and Xiyu Liu}{Positive solutions of singular boundary value problem}

\author{YUXIA WANG and XIYU LIU}

\address{Department of Statistics, Shandong Economic University,
Jinan, Shandong 250 014, People's Republic of China\\
\noindent E-mail: yuxiawty@yahoo.com.cn}

\volume{113}

\mon{May}

\parts{2}

\Date{MS received 29 October 2001; revised 2 January 2003}

\begin{abstract}
This paper investigates the existence of positive solutions of a
singular boundary value problem with negative exponent similar to
standard Emden--Fowler equation. A necessary and sufficient condition for
the existence of $C[0,1]$ positive solutions as well as $C^{1}[0,1]$
positive solutions is given by means of the method of lower and upper
solutions with the Schauder fixed point theorem.
\end{abstract}


\maketitle

\section{Introduction}
Consider the singular boundary value problems for the Emden--Fowler equation
\begin{align}
&u''+p(t)u^{-\lambda}(t) =0,\quad 0 < t < 1,\\[.5pc]
&\alpha u(0)-\beta u'(0)=0,\quad \gamma u(1)+\delta u'(1)=0,
\end{align}
where $ \alpha, \beta, \gamma, \delta\geq 0, \lambda\in R $ and
$\rho:=\gamma\beta +\alpha\gamma+\alpha\delta>0; p\in C((0,1),
[0,\infty))$ and may be singular at $t=0, t=1$. When $\lambda<0$,
see \cite{H,L.P,W,Z.Y} for the result concerning the above problem. When
$\lambda>0$, \cite{T} shows the existence and uniqueness to
(1) and (2) in the case of $\beta=\delta=0$ by means of the
shooting method. For the following problem
\begin{align}
&u^{\prime\prime}+p(t)u^{-\lambda}(t)+q(t)u^{-m}(t)=0, \quad  0<t<1,\\[.5pc]
&\alpha u(0)-\beta u'(0)=0,\quad \gamma u(1)+\delta u'(1)=0,
\end{align}
where $\alpha, \beta, \gamma, \delta \geq 0$,
$\rho=\gamma\beta+\alpha \gamma+\alpha\delta > 0$ and  $p,q
\in((0,1),[0,\infty))$. Mao \cite{M} gave a sufficient and necessary
condition when $\lambda<0,m<0$. In this paper we shall consider the
case of $\lambda>0,m>0$ for the problems (3) and (4).

A function $u(t)\in C^1 [0,1]\cap C^2 (0,1)$ is a positive solution of (3) and (4) if
$u$ satisfies (3) and (4) and $u(t)>0, t\in (0,1). $

\section{Main results}

We state the following hypothesis, which is used throughout
this paper.
\begin{align*}
\hskip -1.7cm\hbox{(H)} \qquad \  p(t), q(t)\in C(0,1), p(t)\!\geq\!0, p(t)& \not\equiv 0,q (t)\!\geq\!0, q(t)
\not\equiv
0,t\in (0,1), \lambda, m>0.
\end{align*}

We now state the main results of this paper as follows:
\begin{theor}[\!]
\ Suppose that {\rm (H)} is satisfied. Then
\end{theor}\vspace{-1pc}
{\it
\begin{enumerate}\leftskip .8pc
\renewcommand{\labelenumi}{\rm (\Roman{enumi})}
\item If $\beta\delta\not=0,$ the problems $(3)$ and $(4)$ have a positive
solution if and only if
\begin{equation}
\hskip -1pc 0 < \int^1_0 [p(t)+q(t)]\,\hbox{\rm d}t < \infty.
\end{equation}
\item  If $\beta=0, \delta\not=0,$ the problems $(3)$ and $(4)$ have a
positive solution if and only if
\begin{equation}
\hskip -1pc 0 < \int^1_0 t[p(t)+q(t)]\,\hbox{\rm d}t < \infty.
\end{equation}
\item If $\beta\not=0, \delta=0,$ the problems $(3)$ and $(4)$ have a
positive solution if and only if
\begin{equation}
\hskip -1pc 0 < \int^1_0 (1-t)[p(t)+q(t)]\,\hbox{\rm d}t < \infty.
\end{equation}
\item If $\beta=\delta=0,$ the problems $(3)$ and $(4)$ have a positive
solution if and only if
\begin{equation}
\hskip -1pc 0 < \int^1_0 t(1-t)[p(t)+q(t)]\,\hbox{\rm d}t < \infty.
\end{equation}
\end{enumerate}}

\begin{theor}[\!]
Suppose that {\rm (H)} is satisfied. Then problems $(3)$ and $(4)$ have a $C^1 [0,1]$ positive
solution if and only if the following inequalities hold.\vspace{.1cm}

\noindent{\rm (H1)}
\begin{equation}
0 < \int^1_0 [t^{-\lambda}p(t)+t^{-m}(t)q(t)]\, \hbox{\rm d}t < \infty, \ \ \beta=0,
\delta\not=0.
\end{equation}
\noindent{\rm (H2)}
\begin{equation}
0 < \int^1_0 [(1-t)^{-\lambda}p(t)+(1-t)^{-m}q(t)]\, \hbox{\rm d}t < \infty,\ \
\beta\not=0, \delta=0.
\end{equation}
\noindent{\rm (H3)}
\begin{equation}
0 <\! \int^1_0 [t^{-\lambda}(1\!-\!t)^{-\lambda}p(t)+t^{-m}(1\!-\!t)^{-m}q(t)]\, \hbox{\rm d}t <
 \infty,\ \beta\!=\!\delta\!=\!0.
\end{equation}
\end{theor}

\section{Proofs of the main results}

First we prove Theorem~2.1. We will prove all the necessary conditions
first then all the sufficient conditions.

\subsection*{1. \it Necessity}

\noindent{\it Case I}: $\beta\delta\not =0.$
 Let $ u(t)\in C^1 [0,1]\cap C^2 (0,1) $ is
a positive solution of (3) and (4). From (4) and the nontrivial
concave function $u(t)$, we know that $u(t)$ must satisfy the
following case:
\begin{equation*}
 u(0) \geq 0,\quad u(1) \geq 0,\quad u'(0) \geq  0,\quad u'(1) \leq 0.
\end{equation*}
Then there exists $ t_0\in [0,1]$ with $u'(t_0)= 0,$
$u^{\prime\prime}(t) < 0$ yield $u'(t)\leq 0, t\in [t_0,1); u'(t)\geq 0,
t\in (0,t_0].$ Let $C_0$ be a constant which satisfies $C_0u(t)<1/4,
t\in[0,1],$ and $1/C_0\geq 4.$ Then 
\begin{align}
p(t)u^{-\lambda}(t) &\geq p(t)(4C_0)^{\lambda},\\[.5pc]
 q(t)u^{-m}(t) &\geq q(t)(4C_0)^m.
\end{align}

By means of (12) and (13), we have
\begin{align*}
u'(t)	&=   \int_t^{t_0} [p(s)u^{-\lambda}(s)+q(s)u^{-m}(s)]\,\hbox{\rm d}s\\[.5pc]
	&\geq \int_t^{t_0} [(4C_0)^{\lambda}p(s)+(4C_0)^m q(s)]\,\hbox{\rm d}s\\[.5pc]
	&\geq (4C_0)^{\min\{\lambda,m\}}\int_t^{t_0} [p(s)+q(s)]\,\hbox{\rm d}s, \ \ t\in[0,t_0)\\[.5pc]
-u'(t) 	&= \int^t_{t_0} [p(s)u^{-\lambda}(s)+q(s)u^{-m}(s)]\,\hbox{\rm d}s \\[.5pc]
	&\geq (4C_0)^{\min\{\lambda,m\}}\int^t_{t_0} [p(s)+q(s)]\,\hbox{\rm d}s,\ \ \  t\in [t_0,1].
\end{align*}
So,
\begin{align*}
  0& < \int^1_0 [p(s)+q(s)]\,\hbox{\rm d}s \\[.5pc]
   & = \int^{t_0}_0 [p(s)+q(s)]\,\hbox{\rm d}s + \int^1_{t_0} [p(s)+q(s)]\,\hbox{\rm d}s\\[.5pc]
   &\leq (4C_0)^{-\min\{\lambda,m\}}[u'(0)-u'(1)]\\[.5pc]
   &<\infty.
\end{align*}

Therefore, (5) holds.\vspace{.3cm}

\noindent{\it Case II}: $\beta = 0, \delta > 0.$ Let $u \in
C[0,1]\cap C^1 (0,1]\cap C^2 (0,1)$ be a positive solution of (3)
and (4). From (4) we obtain $u(0) = 0, u(1) \geq 0, u^{\prime}(1) =
-\gamma\delta^{-1}u(1) \leq 0.$ Then by the concavity of $u$ there
exists $t_0 \in (0,1]$ with $ u^{\prime} (t_0) = 0.$

Let $C_1$ be a constant satisfying $C_1 u(t) \leq {1}/{4}, 
{1}/{C_1} \geq 4.$ By means of (12) and (13), we obtain,
\begin{align*}
 0 &\leq \int^{t_0}_0 s[p(s)+q(s)]\,\hbox{\rm d}s = \int_0^{t_0} \,\hbox{\rm d}\tau\int^{t_0}_{\tau}[p(s)+q(s)]\,\hbox{\rm d}s\\[.5pc]
   &\leq (4C_1)^{-\lambda}\int^{t_0}_0 \,\hbox{\rm d}\tau\int^{t_0}_{\tau} p(s)u^{-\lambda}(s)\,\hbox{\rm d}s\\[.5pc]
&\quad \  +
   (4C_1)^{-m}\int^{t_0}_0 \,\hbox{\rm d}\tau \int^{t_0}_{\tau}q(s)u^{-m}(s)\,\hbox{\rm d}s
\end{align*}
\begin{align*}
   &\leq (4C_1)^{-\min\{\lambda,m\}}\int^{t_0}_0 \,\hbox{\rm d}\tau \int^{t_0}_{\tau}[p(s)u^{-\lambda}
   (s) + q(s)u^{-m}(s)]\,\hbox{\rm d}s\\[.5pc]
   &= (4C_1)^{-\min\{\lambda,m\} }\int^{t_0}_0 (-u^{\prime}(t_0)+u^{\prime}(\tau))\,\hbox{\rm d}\tau\\[.5pc]
   &= (4C_1)^{-\min\{\lambda,m\}}u(t_0)\\[.5pc]
   &< +\infty.
\end{align*}
Similarly,
\begin{equation*}
0 \leq \int^1_{t_0} s[p(s)+q(s)]\,\hbox{\rm d}s < +\infty.
\end{equation*}
Hence we conclude
\begin{equation*}
0 < \int^1_0 t[p(t)+q(t)]\,\hbox{\rm d}t < \infty.
\end{equation*}

\noindent {\it Case III}: $\beta > 0, \delta=0$.
The proof for Case III is almost the same as that for Case II.\vspace{.3cm}

\noindent {\it Case IV}: $\beta = \delta = 0.$ Let $u \in C[0,1]$ be a
positive solution of (3) and (4). Integrating (3) twice gives
\begin{align}
u'\left(\frac{1}{2}\right) - u'(t) = \int^t_{{1}/{2}} [p(s)u^{-\lambda}(s&)+q(s)u^{-m}(s)]\,\hbox{\rm d}s,\nonumber\\[.5pc]
u'\left(\frac{1}{2}\right)\left(t-\frac{1}{2}\right) - u(t) + u\left(\frac{1}{2}\right)& = \int^t_{{1}/{2}}
\,\hbox{\rm d}\eta \int^{\eta}_{{1}/{2}}
[p(s)u^{-\lambda}(s)\nonumber\\[.5pc]
&\quad\, +q(s)u^{-m}(s)]\,\hbox{\rm d}s\nonumber\\[.5pc]
& = \int^t_{{1}/{2}} (t-s)[p(s)u^{-\lambda}(s)\nonumber\\[.5pc]
&\quad \, +q(s)u^{-m}(s)]\,\hbox{\rm d}s.
\end{align}
Since the limit of (23) as $t \arrow 1$ exists and is finite, by the
monotone convergence theorem,
\begin{equation*}
0 < \int^1_{{1}/{2}} (1-s)[p(s)u^{-\lambda}(s)+q(s)
  u^{-m}(s)]\,\hbox{\rm d}s < \infty.
\end{equation*}
So,
\begin{equation*}
0 < \int^1_{{1}/{2}} (1-s)[p(s)+q(s)]\,\hbox{\rm d}s < \infty.
\end{equation*}
Similarly
\begin{equation*}
0 < \int^{{1}/{2}}_0 (1-s)[p(s)+q(s)]\,\hbox{\rm d}s < \infty.
\end{equation*}
Hence
\begin{equation*}
0 < \int^1_0 s(1-s)[p(s)+q(s)]\,\hbox{\rm d}s < \infty.
\end{equation*}

\subsection*{2. \it Sufficiency}

\noindent {\it Case I}: $\beta\delta\not =0.$ Suppose that (5) is
satisfied. Let
\begin{align*}
q_1 (t)&= \frac{\gamma(1-t)+\delta}{\rho} \int^t_0 (\alpha s+\beta)[p(s)+q(s)]
\,\hbox{\rm  d}s\\[.5pc]
&\quad  +\frac{\alpha t+\beta}{\rho}\int^1_t
(\gamma(1-s)+\delta)[p(s)+q(s)]\,\hbox{\rm  d}s,\ t \in[0,1].
\end{align*}
Then $ q_1 \in C^1 [0,1] \cap C^2 (0,1) $ satisfies (4) and solve the
equation $q_1^{\prime\prime}(t)= -$ $[p(t)+q(t)], t\in (0,1).$ Let $L_1=
({\beta\delta}/{\rho}) \int^1_0 [p(s)+q(s)]\,\hbox{\rm  d}s, L_2= 1/\rho \int^1_0
(\alpha s+\beta)$ $(\gamma (1-s)+\delta)[p(s)+q(s)]\,\hbox{\rm  d}s.$ Then it is easy to
check that $0\!<\!L_1\!\leq\!q_1 (t) \leq L_2,$ $ \forall t\in [0,1].$ Let
$\alpha(t)\!\!=\!\!k_1q_1 (t),$ $\beta(t)\!\!=\!\!k_2q_1 (t), t\in [0,1],$ where $k_1\!\!=\!\!
\min$ $\{1, L_2^{{-\lambda}/{(1+\lambda)}}, L_2^{{-m}/{(1+m)}}\}, k_2=
\max\{1, L_1^{{-\lambda}/{(1+\lambda)}}, L_1^{{-m}/{(1+m)}}\}. $Then
$\alpha(t), \beta(t) \in C^1$ $[0,1] \cap C^2 (0,1), 0 < \alpha(t) \leq
\beta(t),$ $t\in [0,1], $ and $\alpha(t), \beta(t)$ satisfy the boundary
condition (4). Furthermore,
\begin{align*}
\alpha^{{\prime\prime}}(t)+ p(t)\alpha^{-\lambda}(t)+ q(t)\alpha^{-m}(t)
&= -k_1[p(t)+q(t)]+ p(t)[k_1 q_1 (t)]^{-\lambda}\\[.5pc]
&\quad \ + q(t)[k_1 q_1 (t)]^{-m} \\[.5pc]
&\quad \ \geq p(t)[(k_1L_2)^{-\lambda}- k_1]\\[.5pc]
&\quad \ + q(t)[(k_1L_2)^{-m} - k_1]\\[.5pc]
&\quad \ \geq 0,\ \ t\in(0,1),\\[.5pc]
\beta^{{\prime\prime}} (t)+ p(t)\beta^{-\lambda}(t)
+ q(t)\beta^{-m}(t)
&=- k_2[p(t)+q(t)]+ p(t)[k_2 q_1(t)]^{-\lambda}\\[.5pc]
&\quad \ + q (t)[k_2q_1(t)]^{-m}\\[.5pc]
&\quad \ \leq p(t)[(k_2L_1)^{-\lambda}-k_2]\\[.5pc]
&\quad \ + q(t)[(k_2L_1)^{-m}-k_2]\\[.5pc]
&\quad \ \leq 0,\ \ t\in(0,1).
\end{align*}
Thus, $\alpha(t)$ and $\beta(t)$ are respectively lower and upper solutions of
 problems (3) and (4).
We will now prove that problems (3) and (4) admit a $C^1 [0,1]$ positive
solution $u^{*}$ satisfying $0< \alpha(t)\leq u^{*}(t)\leq \beta(t), t\in [0,1].$

First, define an auxiliary function
\begin{equation*}
f(t,u)=\begin{cases}
p(t)\alpha^{-\lambda}(t)+ q(t)\alpha^{-m}(t), & u < \alpha(t), \\[.5pc]
p(t)u^{-\lambda}(t)+ q(t)u^{-m}(t), & \alpha\leq u\leq \beta(t),\\[.5pc]
p(t)\beta^{-\lambda}(t)+ q(t)\beta^{-m}(t), & u> \beta(t).
\end{cases}
\end{equation*}
From (H), $f: (0,1)\times R\arrow [0,\infty)$ is continuous.

Consider the boundary value problem
\begin{align}
&u^{\prime\prime} (t)+ f(t,u(t))= 0,\quad t\in (0,1),\\[.5pc]
&\alpha u(0)- \beta u^{\prime}(0)= 0,\quad  \gamma u(1)+ \delta u^{\prime}(1)= 0.
\end{align}
It is clear that the above problem is equivalent to the integral equation
\begin{equation}
u(t)= Au(t)= \int^1_0 G(t,s)f(s,u(s))\,\hbox{\rm d}s,
\end{equation}
where
\begin{equation*}
G(t,s)=\frac{1}{\rho}\begin{cases}
(\alpha s+\beta)(\gamma(1-t)+\delta), \ s< t,\\[.5pc]
(\alpha t+\beta)(\gamma(1-s)+\delta), \ t\geq s.
\end{cases}
\end{equation*}
Let $ X = C[0,1].$ For $u \in X,$ if for some $t\in [0,1]\  \alpha(t)\leq u(t)
\leq \beta (t)$, we obtain that
\begin{align}
0&\leq p(t)u^{-\lambda}(t)+ q(t)u^{-m}(t) \leq p(t)\alpha^{-\lambda}(t)+ q(t)
\alpha^{-m}(t)\nonumber\\[.5pc]
&\leq p(t)(k_1L_1)^{-\lambda}+ q(t)(k_1L_1)^{-m}\leq (k_1L_1)^{-\min\{\lambda,m\}}
[p(t)+q(t)].
\end{align}
Therefore, from (5), (14) and (18), we know that $ A: X \arrow X$ is
continuous and $A(X)$ is a bounded set. In addition, $u\in X \cap C^1
[0,1]$ is a solution of problems (15) and (16) if and only if $ Au=u $.

Since $AX \subseteq C^2[0,1],$ by the standard application of Arzela--Ascoli
theorem, we obtain that $A$ is compact. By means of Schauder fixed point
theorem, we obtain that $A$ has at least one fixed point $u^{*}\in
X\cap C^1 [0,1].$ We will show
\begin{equation}
\alpha(t)\leq u^{*}(t)\leq \beta(t),\quad  t\in[0,1],
\end{equation}
which will imply that $u^{*}(t)\in C^1 [0,1]$ is a positive solution of
(3) and (4). Suppose (19) is not satisfied. Then there exists $ t^{*} \in
[0,1]$ such that either $u^{*}(t^*) < \alpha(t^{*})$ or $ u^{*}(t^*) >
\beta(t^{*}). $ Let us consider the second case. Let $I \subseteq [0,1]$
denote the maximal interval containing $t^{*}$ such that $u > \beta$ on
I. Then, it is clear that either $u=\beta$ on $\partial I$ or both $u$
and $\beta$ satisfy the same boundary conditions given by (4) on
$\partial I$. Let $Z(t) = \beta(t)-u^{*}(t), t \in I.$ Then
$Z^{\prime\prime} (t) \leq 0$ on $I$ and either $Z(t) = 0$ for $t\in
\partial I$ or $Z(t)$ satisfies the boundary condition (4) for $t\in
\partial I$. From the maximum principle $Z(t) \geq 0$ for $t \in I, \hbox{i.e.}\
\beta(t) \geq u^{*}(t)$ for  $t \in I.$ This is a contradiction. In the same
way, $\alpha (t) \leq u^{*}(t)$ for $t\in I.$ So, $u^{*}(t)$ is a $C^1 [0,1]$
positive solution of (3) and (4).\vspace{.3cm}

\noindent {\it Case II}: $\beta = 0, \delta > 0.$ \ Suppose (6) is
satisfied. Choose $n \geq 4$ so that $n\min\{\lambda,m\} > 1$. Let
\begin{align*}
R(t) &= \left(\frac{\gamma (1-t)+\delta}{\gamma+\delta}\int^t_0 s[p(s)+q(s)]\,\hbox{\rm d}s\right.\\[.5pc]
&\left.\quad \ +\  t \int^1_t \frac{\gamma
(1-s)+\delta}{\gamma+\delta}[p(s)+q(s)]\,\hbox{\rm d}s\right)^{1/(n\min\{\lambda,m\})}.
\end{align*}
Then $ R(t) \in C[0,1] \cap C^2 (0,1),$ satisfies $R(t) > 0,
R^{\prime\prime} (t) \leq 0, t \in (0,1), $ and the boundary conditions
$R(0)=0, R(1) > 0, \gamma R(1)+\delta R^{\prime}(1) \geq 0$. We now
estimate
\begin{align*}
0 &\leq \int^1_0 t[p(t)+q(t)]R^{-\min\{\lambda,m\}}(t)\,\hbox{\rm d}t\\[.5pc]
  &= \int^1_0 t[p(t)+q(t)]\Bigg(\frac{\gamma(1-t)+\delta}{\gamma+\delta}\int^t_0
s[p(s)+q(s)]\,\hbox{\rm d}s\\[.5pc]
&\quad \  + t\int^1_t \frac{\gamma(1-s)+\delta}{\gamma+\delta}[p(s)+q(s)]
\,\hbox{\rm d}s\Bigg)^{-({1}/{n})}\,\hbox{\rm d}t\\[.5pc]
&\leq \int^1_0 t[p(t)+q(t)] \left( \int^1_0
s[p(s)+q(s)]\,\hbox{\rm d}s\right)^{-({1}/{n})}\,\hbox{\rm d}t.
\end{align*}
Let 
\begin{align*}
\Gamma_1 (t) &=
\frac{\gamma(1-t)+\delta}{\gamma+\delta}\int^t_0 s[p(s)+q(s)] \,\hbox{\rm d}s\\[.5pc]
&\quad\ +
t\int^1_t \frac{\gamma(1-s)+\delta}{\gamma+\delta}[p(s)+q(s)]\,\hbox{\rm d}s,\\[.5pc]
\Gamma_2(t) &= \frac{\gamma(1-t)+\delta}{\gamma+\delta}\int^t_0 s[p(s)+q(s)]
R^{-\min\{\lambda,m\}}\,\hbox{\rm d}s\\[.5pc]
&\quad \  + t\int^1_t \frac{\gamma(1-s)+\delta}{\gamma+\delta}
[p(s)+q(s)]R^{-\min\{\lambda,m\}}(s)\,\hbox{\rm d}s\\[.5pc]
&\quad \ +R(t),\ \ t \in [0,1].
\end{align*}
It is clear that $\gamma\Gamma_1(1)+\delta\Gamma_1^{\prime}(1) =
0,\gamma\Gamma_2(1) +\delta\Gamma_2^{\prime}(1) \geq 0$ and $L_3
t({\gamma(1-t)+\delta}/{(\gamma +\delta)}) \leq \Gamma_1 (t) \leq L_3,
R(t) \leq \Gamma_2 (t) \leq L_4, t\in [0,1],$ where
\begin{align*}
L_3 &= \int^1_0 s[p(s)+q(s)]\left(\frac{\gamma(1-s)+\delta}{\gamma+\delta}\right)\,\hbox{\rm d}s,\\[.5pc]
L_4 &= \int^1_0 s\left(\frac{\gamma(1-s)+\delta}{\gamma+\delta}\right)[p(s)+q(s)]R^{-\min
\{\lambda,m\}}(s)\,\hbox{\rm d}s\\[.5pc]
&\quad \ + R_0, \ R_0 = \max_{t\in [0,1]}R(t).
\end{align*}
We also check by direct computation that $\Gamma_1^{\prime\prime}(t) =
-[p(t)+q(t)],\ \Gamma_2^{\prime\prime}(t) \leq
-(p(t)+q(t))R^{-\min\{\lambda,m\}}(t),\ \ t \in (0,1)$. Let $\alpha (t)
= k_1 \Gamma_1 (t), \beta(t) = k_2 \Gamma_2 (t), t \in [0,1],$ where $
k_1 = \min\{1, L_3^{-\lambda/1+\lambda}, L_3^{-m/1+m}\}, k_2 = \max\{1,
L_4^{\min \{\lambda,m\}}\}.$ Then we have
\begin{align*}
\alpha^{\prime\prime}(t)+p(t)\alpha^{-\lambda}(t)+q(t)\alpha^{-m}(t)
&\geq -k_1 [p(t)+q(t)]+p(t)(k_1L_3)^{-\lambda}\\
&\quad \ +q(t)(k_1L_3)^{-m} \geq 0,
\end{align*}
\begin{align*}
\beta^{\prime\prime}(t)+p(t)\beta^{-\lambda}(t)
  +q(t)\beta^{-m}(t)
&\leq -k_2[p(t)+q(t)]R^{-\min\{\lambda,m\}}(t)\\
&\quad \ +p(t)(k_2
R(t))^{-\lambda}(t)\\
&\quad \ +q(t) (k_2R(t))^{-m} \leq 0.
\end{align*}
In addition, $\alpha(t), \beta(t) \in C[0,1] \cap C^1 [0,1] \cap
C^2 (0,1), \gamma\alpha(1)+\delta\alpha^{\prime}(1) = 0, \gamma\beta(1)
+ \delta\beta^{ \prime}(1) \geq 0.$ Let $Z(t) = \beta(t)-\alpha(t).$
Then $Z(0) = 0, \gamma Z(1) + \delta Z^{\prime}(1) \geq 0.$ Also, if we
assume $ \alpha \leq \beta$ on an interval $I \subseteq [0,1],$ we find
$Z^{\prime\prime}(t) \leq -k_2[p(t)+q(t)]R^{-\min\{\lambda,m\}} +
k_1[p(t)+q(t)] \leq -k_2 [p(t)+q(t)]L_4^{-\min\{\lambda,m\}} +
[p(t)+q(t)] \leq 0$ in $I$. As before, for $t \in \partial I$, we have
either $ Z(t) = 0$ or if $t=1, \gamma Z(1) + \delta Z^{\prime}(1) \geq
0.$ From the maximum principle, $Z(t) \geq 0, t \in [0,1],$ i.e.
$\alpha(t) \leq \beta(t), t \in [0,1],$ which is a contradiction. Hence
$\alpha(t), \beta(t)$ are respectively the lower and upper solutions of
(3) and (4).

In the following we will prove that problems (3) and (4) have a $C^1
[0,1]$ positive solution $u(t)$ satisfying  $0 < \alpha(t) \leq
u(t) \leq \beta(t), t\in(0,1).$ Let ${a_n}$ be a sequence
satisfying $0 < \cdots <a_{n+1} < a_n < \cdots < a_2 < a_1 < 1/2$
with $ a_n \arrow 0 $ as $ n \arrow \infty.$ Let ${r_n}$ be a
sequence satisfying $ 0 < \alpha(a_n) \leq r_n \leq \beta(a_n),
n=1,2,\ldots$\,. For each $n$, consider the following singular
boundary problem:
\begin{align}
&u^{\prime\prime}(t) + p(t)u^{-\lambda}(t) + q(t)u^{-m}(t) = 0,\quad t \in (a_n,1),\\[.5pc]
&u(a_n) = r_n,\quad \gamma u(1) + \delta u^{\prime}(1) = 0.
\end{align}
From (6), we see that
\begin{equation*}
\int^1_0 s[p(s)+q(s)]\,\hbox{\rm d}s \geq \int^1_{a_n} s[p(s)+q(s)]\,\hbox{\rm d}s \geq a_n
\int^1_{a_n} [p(s)+q(s)]\,\hbox{\rm d}s.
\end{equation*}
Therefore,
\begin{equation}
0 \leq \int^1_{a_n} [p(s)+q(s)]\,\hbox{\rm d}s \leq \frac{1}{a_n}\int^1_0 s[p(s)+q(s)]\,\hbox{\rm d}s <
\infty.
\end{equation}
Following the proof of Case I, with (22), we can say that for each $n$,
the singular boundary value problems (20) and (21) have at least one positive
solution $u_n(t) \in C^1 [a_n,1] $ satisfying $\alpha(t) \leq u_n (t)
\leq \beta(t), t \in [a_n,1]$. So
\begin{align*}
&|u_n(1)| \leq M = \alpha(1) + \beta(1),\\[.5pc]
&|u_n^{\prime}(1)| = |(\gamma/\delta)u_n(1)| \leq \frac{\gamma}{\delta}M,\ \  n = 1,
2,\ldots.
\end{align*}
Without loss of generality, we can assume
\begin{align*}
u_n(1) &\arrow u_0 \in [\alpha(1),\beta(1)], \ \  n \arrow \infty,\\[.5pc]
u_n^{\prime}(1) &\arrow -(\gamma/\delta) u_0, \ \ n \arrow \infty.
\end{align*}
Similar to Theorem~3.2 in \cite{H}, we can prove (3) has a positive solution
$u(t)$ with $u(1)=u_0, u^{\prime}(1)=-(\gamma/\delta) u_0.$ Its maximal
interval of existence is $(\omega^-,\omega^+)$, and $u_n(t)$ converges
to $ u(t)$ uniformly in any compact subset of $(\omega^-,\omega^+)$
($u_n^{\prime}(t) $ converges to $u^{\prime}(t)$ uniformly), $n \arrow
\infty$. Since $\alpha(t) \leq u_n (t) \leq \beta (t), t \in [a_n,1]$
and $\bigcup^{\infty}_{n=1} [a_n,1] = [0,1],$ we have $\alpha(t) \leq
u(t) \leq \beta(t), t \in (\omega^-,\omega^+) \cap [0,1].$ From
continuation theorem we obtain that $[0,1] \subset (\omega^-,\omega^+)$.
Since $\alpha(0) = \beta(0) = 0, $ we also obtain that $u(0) = 0$. In
addition $\delta u(1) + \gamma u^{\prime}(1) = 0$. Thus $u(t) \in C[0,1]
\cap C^1 (0,1) \cap C^2 (0,1)$ is a positive solution of 
problems (3) and (4). The proof for Case II is complete.\vspace{.3cm}

\noindent {\it Case III}: $\beta > 0, \delta = 0.$ The proof for Case
III is almost the same as that for Case II.\vspace{.3cm}

\noindent {\it Case IV}: $\beta = \delta = 0.$ Let $ Q(t) =
((1-t)\int^t_0 s[p(s)+q(s)]\,\hbox{\rm d}s + t\int^1_t (1-s)
[p(s)$\break $+q(s)]\,\hbox{\rm d}s)^{{1}/{n\min\{\lambda,m\}}},\ \ t\in [0,1].$ Then
$Q(t) \in C[0,1] \cap C^2 (0,1)$ with $Q(t) > 0,$\break $ Q^{\prime\prime}(t)
\leq 0, t \in (0,1)$ and $Q (0) = Q (1) = 0$. Let 
\begin{align*}
\Gamma_1 (t) &= (1-t)\int^t_0 s[p(s)+q(s)]\,\hbox{\rm d}s + t\int^1_t
(1-s)[p(s)+q(s)]\,\hbox{\rm d}s,\\[.5pc]
\Gamma_2 (t) &= (1-t)\int^t_0 s[p(s)+q(s)]Q^{-\min\{\lambda,m\}}(s)\,\hbox{\rm d}s\\[.5pc]
&\quad \  + t
\int^1_t (1-s)[p(s)+q(s)]Q^{-\min\{\lambda,m\}}(s)\,\hbox{\rm d}s + Q(t) \ t
\in [0,1].
\end{align*}
Then $t(1-t)L_5 \leq \Gamma_1 (t) \leq L_5, Q(t) \leq \Gamma_2(t) \leq
L_6, L_5 = \int^1_0 s(1-s)[p(s)+q(s)]\,\hbox{\rm d}s, \ L_6 = \int^1_0
s(1-s)[p(s)+q(s)] Q^{-\min\{\lambda,m\}}(s)\,\hbox{\rm d}s + Q_0,\ \ Q_0 = \max Q(t).
$ Let $\alpha(t) = k_1 \Gamma_1 (t), \beta(t) = k_2 \Gamma_2 (t), t \in
[0,1]$. Here
\begin{equation*}
k_1 = \min\{ 1, L_5^{-\lambda/1+\lambda}, L_5^{-m/1+m}
\},\quad k_2 = \{ 1, L_6^{\min \{\lambda,m\}} \}.
\end{equation*}
Then $\alpha(t), \beta(t)$ are respectively the lower and upper solutions of
(3) and (4).The remaining proof is analogous to that of Case II. Thus
we complete the proof of Theorem~2.1.\vspace{.3cm}

\noindent{\it Proof of Theorem}~2.2\vspace{.3cm}

\noindent Proof for Case (H1): $\beta = 0,  \delta > 0$.

\subsection*{1. {\it Necessity}}

Suppose that $u$ is a $C^1 [0,1]$ positive solution of (3) and (4). Then
both $u'(0)$ and $u'(1)$ exist, and $p(t), q(t) \not\equiv 0, t \in
(0,1).$ From (4) and the fact that $u$ is a positive concave function,
we know that $u(0) = 0, u(1) > 0, u'(0) > 0, u'(1) \leq 0.$ Then there
exists $t_0 \in (0,1]$ such that $u'(t_0) = 0$. Since
$u^{\prime\prime}(t) \leq 0,$ and $u(t) > 0,\ t \in [0,1]$, one easily sees
that $0 < u(1) \leq u(t) \leq u(t_0),\ t \in [t_0,1]$ and there
exist constants $I_1$ and $I_2$ which satisfy
\begin{equation}
I_1 t \leq u(t) \leq I_2 t,\ \ t \in [0,1].
\end{equation}
Hence 
\begin{align*}
\hskip -1pc 0 &< \int^1_0 [p(s)s^{-\lambda}+q(s)s^{-m}]\,\hbox{\rm d}s \leq I_2^{\max\{\lambda,m\}} \int^1_0 [p(s)u^{-\lambda}(s)+q(s)u^{-m}(s)]\,\hbox{\rm d}s
\end{align*}
\begin{align*}
&= -I_2^{\max\{\lambda,m\}}\int^1_0 u^{\prime\prime}(s)\,\hbox{\rm d}s
= I_2^{\max\{\lambda,m\}}(u'(0)-u'(1)) < \infty.
\end{align*}
The above inequality shows that (9) holds.

\subsection*{2. {\it Sufficiency}} 

Suppose that (9) is satisfied. Let
\begin{align*}
\Gamma (t) &=
\frac{\gamma(1-t)+\delta}{\gamma+\delta}\int^t_0 s[p(s)s^{-\lambda}
+q(s)s^{-m}]\,\hbox{\rm d}s\\[1pc]
&\quad \ + t\int^1_t
\frac{\gamma(1-s)+\delta}{\gamma+\delta}[p(s)s^{-\lambda}
+q(s)s^{-m}]\,\hbox{\rm d}s, t \in [0,1].
\end{align*}
Then $\Gamma(t) \in C^1 [0,1] \cap C^2 (0,1)$. Replace $u(t)$ with
$\Gamma(t)$ in (24) and let
\begin{align*}
I_1 &= \frac{\delta}{\gamma+\delta}\int^1_0
s\left(\frac{\gamma(1-s)+\delta}{\gamma+\delta}\right) [p(s)s^{-\lambda}
+q(s)s^{-m}]\,\hbox{\rm d}s,\\[.8pc]
I_2 &= \int^1_0 [p(s)s^{-\lambda}+q(s)s^{-m}]\,\hbox{\rm d}s.
\end{align*}
Then $\Gamma(t)$ satisfies (24). Let
\begin{align*}
&k_1 = \min\{ 1, (I_2)^{-\lambda/1+\lambda}, (I_2)^{-m/1+m} \},\\  
&k_2 = \max\{ 1, (I_1)^{-\lambda/1+\lambda}, (I_1)^{-m/1+m} \}.
\end{align*}
Then $\alpha(t) = k_1\Gamma(t), \beta(t) = k_2\Gamma(t),\ t \in [0,1].$
 Then $\alpha(t)$ and $\beta(t)$ are respectively lower and upper solutions of problems (3) and (4).
 It is clear that $0 < \alpha(t) \leq \beta(t), t \in (0,1], \alpha(0) = \beta(0) = 0,
 \gamma \alpha(1)+\delta\alpha'(1) = 0, \gamma\beta(1)+\delta\beta'(1) = 0.$
 On the other hand, when $t \in (0,1), \alpha(t) \leq u \leq \beta(t)$, we have
\begin{align*}
0 \leq p(t)u^{-\lambda}(t)+q(t)u^{-m}(t) &\leq p(t)\alpha^{-\lambda}(t)+
 q(t)\alpha^{-m}(t)\\[.8pc]
&= p(t)(k_1 I_1 t)^{-\lambda} + q(t)(k_1 I_1 t)^{-m} = F(t).
\end{align*}
From (9), we have $\int^1_0 F(t)\,\hbox{\rm d}t < \infty$. The same argument that
we have given in the sufficiency of Theorem 2.1 assures us that problems
(3) and (4) admit a positive solution $u \in C^1 [0,1] \cap C^2 (0,1)$
such that $\alpha(t) \leq u(t) \leq \beta(t),\ t \in [0,1]$. The proof
for Case (H1) is complete. Similarly we can prove cases 
(H2) and (H3). Thus we complete the proof of Theorem 2.2.

\section*{Acknowledgement}

This work is supported in part by the NSF(Youth) of Shandong Province
 and NNSF of China.

\end{document}